\documentclass{elsart}
\usepackage[utf8]{inputenc}
\usepackage{amsmath}
\usepackage{amsfonts}
\usepackage{amssymb}
\usepackage{amstext}
\usepackage{amsgen}
\usepackage{amsbsy}
\usepackage{amsopn}
\newtheorem{theorem}{Theorem}[section]
\newtheorem{lemma}[theorem]{Lemma}
\newtheorem{proposition}[theorem]{Proposition}

\newtheorem{definition}[theorem]{Definition}

\newtheorem{conjecture}[theorem]{Conjecture}
\newtheorem{remark}[theorem]{Remark}
\newtheorem{corollary}[theorem]{Corollary}
\newtheorem{example}[theorem]{Example}
\newcommand{\supp}{\mathop{\mathrm{supp}}}

\newenvironment{proof}{\emph{Proof.}} {\qed}
\begin{document}
\def\C{{\mathbb C}}
\def\N{{\mathbb N}}
\def\R{{\mathbb R}}
\def\Z{{\mathbb Z}}
\def\Q{{\mathbb Q}}

\def\abs#1{\left|#1\right|}
\def\Cnt#1{{\mathcal C}^{#1}}
\def\comp{\circ}
\def\conv{\star}
\def\csub{\subset\subset}
\def\defstyle#1{{\bf #1}}
\def\eps{\varepsilon}
\def\Gen{{\mathcal G}}
\def\implies{\Rightarrow}
\def\integr{\int\limits}
\def\kar#1{\chi_{#1}}
\def\Null{{\mathcal N}}
\def\uniform{\rightrightarrows}
\begin{frontmatter}
	\title{Homogeneity in generalized function algebras}
	\author{Clemens Hanel \& Eberhard Mayerhofer}
	\address{Faculty of Mathematics, University of Vienna, Nordbergstra\ss e 15, 1090 Vienna, Austria}
	\ead{clemens.gregor.hanel@univie.ac.at}
	\ead{eberhard.mayerhofer@univie.ac.at}

	\author{Stevan Pilipovi\'c}
	\address{Faculty of Sciences and Mathematics, University of Novi Sad, Trg Dositeja Obradovica 4, 21000 Novi Sad, Serbia}
	\ead{pilipovic@im.ns.ac.yu}
	
	\author{Hans Vernaeve}
	\address{Faculty of Civil Engineering, University of Innsbruck, Technikerstra\ss e 13, 6020 Innsbruck, Austria}
	\ead{hans.vernaeve@uibk.ac.at}
\begin{abstract}
We investigate homogeneity in the special Colombeau algebra on
$\mathbb R^d$ as well as on the pierced space $\mathbb
R^d\setminus\{0\}$. It is shown that strongly scaling invariant
functions on $\mathbb R^d$ are simply the constants. On the pierced
space, strongly homogeneous functions of degree $\alpha$ admit
tempered representatives, whereas on the whole space, such functions
are polynomials with generalized coefficients. We also introduce
weak notions of homogeneity and show that these are consistent with
the classical notion on the distributional level. Moreover, we
investigate the relation between generalized solutions of the Euler
differential equation and homogeneity.
\end{abstract}
\end{frontmatter}

\section{Introduction}\label{chapterconstant}
Differential algebras of generalized functions in the sense of J.\
F.\ Colombeau (cf.\ \cite{Colombeau, C}) have proved valuable as a
tool for treating partial differential equations with singular data
or coefficients. Over the past twenty years a large number of
applications have been published in scientific journals (cf.\
\cite{MObook}). Also, the theory has been adapted to allow for
applications in a geometric context as well as in relativity (cf.\ \cite{Bible} and the references therein). As a natural consequence of intense research
in the field, increasing importance is ascribed to an understanding
of  intrinsic problems in algebras of generalized functions. This is
emphasized by a number of scientific papers on algebraic (cf.\
\cite{A1,algfound}) and topological topics (cf.\ \cite{DHPV,
Scarpalezos1, Garetto2, Garetto1}). 

In this article we characterize 
homogeneous generalized functions in the special Colombeau algebra $\mathcal G$. We investigate this
on the level of coupled calculus (\cite{Bible}, chapter 1), that is we
consider homogeneity (1) in the strong sense, i.e., we solve (H) $u(\lambda x)=\lambda^\alpha u(x)$ in $\mathcal G$ as well as
(2) in the weak sense, which leads to results consistent with distribution theory.

Solving identities as (H) in $\mathcal G$ is non-trivial, for elements of $\mathcal G$ both have
$C^\infty$ character, since they are basically nets of smooth functions 
but also distributional character, since they are not determined pointwise (see section \ref{presec}). In view of
the definition of generalized functions as nets of smooth functions with a prescribed asymptotic growth behaviour with respect
to the smoothing parameter (see section \ref{presec}), a precise asymptotic study
is required for solving equality (H), cf.\ the proof of Theorem \ref{mainthhom}.

We shall see that (1) equality merely yields polynomial solutions (with generalized coefficients, cf.\ Theorem \ref{mainthhom}) as in the $C^\infty$ case,
whereas (2) in the weak sense, all homogeneous distributions are recovered (see for instance Proposition \ref{distshadow}). However,
we further show that there exist generalized functions (weakly) homogeneous of degree $\alpha$ with no distributional shadow
(Proposition \ref{moreweaksolutions}). 

In this context we should note that by embedding distributions into an algebra of Colombeau generalized functions some properties of distributions
are lost in the strong sense. Of course they are not lost in the weak sense, that is on the level of association. This confirms the
statement that Colombeau generalized functions are a natural generalization of Schwarz distributions if we work on the level of association rather than
on the level of equalities (in $\mathcal G$). Thus some equations can have distributional solutions but not the corresponding solutions in $\mathcal G$
(cf.\ \ref{moreweaksolutions}).

It sounds paradoxical but we have ''more'' homogeneous distributions in $\mathcal D'$ than in $\mathcal G$. Our paper is about questions on
strong and associated homogeneity of generalized functions since they are surprisingly different.

\subsection{Program of this paper}
We shall recall a couple of preliminary statements in section \ref{presec}. In particular, we revisit
point value questions and the partial order $\leq$ on the ring of generalized numbers. The aim of section \ref{mainsec}
is to show our main result, that any generalized function which is invariant under standard scaling, is constant. The final section
\ref{homsec} investigates three different concepts of homogeneity in the special Colombeau algebra. We also check consistency of the latter
with the distributional notion.
\section{Preliminaries}\label{presec}
To begin with, we recall introductory material on generalized functions (for more information we refer to \cite{Bible}, chapter 1). Let $\Omega$ be an open subset of $\mathbb R^d$ ($d\geq 1$). In the following we denote by $K\subset\subset \Omega$ a compact set in $\Omega$. The {\it special Algebra} $\mathcal G(\Omega)$ due to J.\ F.\ Colombeau is given by the quotient
\[
\mathcal G(\Omega):=\mathcal E_M(\Omega) /\mathcal N(\Omega),
\]
where the ring of moderate functions $\mathcal E_M(\Omega) $ resp.\ the ring of negligible elements (being an ideal
in $\mathcal E_M(\Omega) $) is given by
\begin{eqnarray}\nonumber
\mathcal E_M(\Omega)&:=&\{ (u_\varepsilon)_\varepsilon\in C^{\infty}(\Omega)^{(0,1]}\, |\, \forall\, K\subset\subset \Omega \,\forall \alpha\in\mathbb N_0^d \, \exists\, N\geq 0 \\\nonumber
& &\sup_{x\in K} |\partial^\alpha u_\varepsilon(x)| = O(\varepsilon^{-N}),\,(\varepsilon\rightarrow 0)\}
\end{eqnarray}
resp.\
\begin{eqnarray}\nonumber
\mathcal N(\Omega)&:=& \{(u_\varepsilon)_\varepsilon\in C^{\infty}(\Omega)^{(0,1]}\, |\, \forall\, K\subset\subset \Omega \,\forall \alpha \in\mathbb N_0^d\, \forall\, m\geq 0\\\nonumber
& &\sup_{x\in K} |\partial^\alpha u_\varepsilon(x)| = O(\varepsilon^{m}),\,(\varepsilon\rightarrow 0)\}.
\end{eqnarray}

The algebraic operations ($+,\,\cdot$) as well as (partial)
differentiation, composition of functions etc. are meant to be
performed component-wise on the level of representatives; the
transfer to the quotient $\mathcal G(\Omega)$ is then well
defined (cf.\ the comprising presentation in the first chapter of
\cite{Bible}). Once a Schwartz mollifier $\rho$ on $\mathbb R^d$
with all moments vanishing has been chosen, the space of compactly
supported distributions may be embedded  canonically into $\mathcal
G(\Omega)$ via convolution; an embedding $\iota_\rho$ of all of $\mathcal
D'(\Omega)$ into our algebra is achieved via a relatively
compact partition of unity using sheaf theoretic arguments. We
remind that this embedding is not canonical since it heavily depends
on the choice of mollifier $\rho$.
Recall that $\iota_\rho$ is a linear embedding which preserves partial differentiation,
multiplication of $C^\infty$ functions and restrictions to open subsets.
\subsection{The ring of generalized numbers}
 Next, we introduce the set of moderate numbers $\widetilde{\Omega}$, defined by the quotient
\[
\widetilde{\Omega}:=\Omega_M/\sim,
\]
where the set of moderate numbers
\[
\Omega_M := \{(x_\varepsilon)_\varepsilon\in \Omega^{(0,1]}\,|\,
\exists\, N\geq 0:\; \Vert x_\varepsilon\Vert = O(\varepsilon^{-N}),\,(\varepsilon\rightarrow 0)\},
\]
 and $\sim$ is an equivalence relation on $\Omega_M$ defined by
\[
(x_\varepsilon)_\varepsilon\sim(y_\varepsilon)_\varepsilon\Leftrightarrow \forall\, p\geq 0:\,\Vert x_\varepsilon-y_\varepsilon\Vert=O(\varepsilon^p),\,(\varepsilon\rightarrow 0).
\]

If $(x_\varepsilon)_\varepsilon\in \Omega_M$, then we denote the class of the latter by $[(x_\varepsilon)_\varepsilon]$. Let $\widetilde{\Omega}_c$ denote the set of compactly supported elements of $\widetilde{\Omega}$, that is: $x_c$ lies
in $\widetilde{\Omega}_c$ if and only if for one (hence any)
representative $(x_\varepsilon)_\varepsilon$ of $x$ there exists a compact
set $K\subseteq \Omega$  and some index
$\varepsilon_0$ such that for all $\varepsilon<\varepsilon_0$ we
have $x_\varepsilon\in K$. Note that $\Omega\subset\widetilde{\Omega}_c\subset \widetilde{\Omega}$, the first inclusion being understood
through the constant embedding $\sigma_0:\Omega\rightarrow\widetilde{\Omega}$, $x\mapsto [(x)_\varepsilon]$.

Let $\mathbb K$ be $\mathbb R$ or $\mathbb C$. If $\Omega=\mathbb K$ in the above definition then we call $\widetilde{\mathbb K}:=\widetilde{\Omega}$ the ring of generalized numbers. Similarly,  $\widetilde{\mathbb K}_c$
is called the ring of compactly supported generalized numbers. Note that $\widetilde{\mathbb C}=\widetilde{\mathbb R}+i\widetilde{\mathbb R}$. For $\Omega=\mathbb R^+$ we set $\widetilde{\mathbb R}^+=\widetilde{\Omega}$, and the compactly supported points in the latter we denote by $(\widetilde{\mathbb R}^+)_c$. We observe that $\widetilde{\Omega}$ and $\widetilde{\Omega}_c$ are rings, whenever $\Omega$  is a subring of $\mathbb R$ and further $\widetilde{\Omega}$ and $\widetilde{\Omega}_c$ are modules over $\widetilde{\mathbb R}$ whenever
$\Omega$ is a vector subspace of $\mathbb R^d$. In particular, $\widetilde{\mathbb R}^d$ is a module over $\widetilde{\mathbb R}$. Finally, note that $\widetilde{\mathbb R^d}=(\widetilde{\mathbb R})^d$.

It can easily be shown that evaluation of generalized functions $f$ on compactly supported generalized points $x_c$
makes perfect sense in the following way: let $(f_\varepsilon)_\varepsilon$ be a representative of
 $f\in\mathcal G(\Omega)$, then
\[
\widetilde f(x_c):=[(f_\varepsilon(x_\varepsilon))_\varepsilon]
\]
yields a well defined generalized number. We denote by $\widetilde
f:\widetilde{\Omega}_c\rightarrow\widetilde{\mathbb C}$ the above map induced by the generalized function $f$.

By a standard point $x$ we shall mean an element of
$\widetilde{\Omega}$ which admits a constant representative, $(\alpha)_\varepsilon$, $\alpha$ being a fixed point in $\Omega$. M.\ Kunzinger and M.\ Oberguggenberger show in
(\cite{MO1}) that it does not suffice to know the values of
generalized functions at standard points in order to determine them
uniquely. 

To see this, take some $\varphi\geq 0\in\mathcal D(\mathbb R)$ with $\supp\varphi\in[-1,1]$ and $\int \varphi=1$ and set $u_\varepsilon:=\varphi_\varepsilon(x-\varepsilon)$,
where $\varphi_\varepsilon(y):=\frac{1}{\varepsilon}\varphi(\frac{y}{\varepsilon})$. Then $(u_\varepsilon)_\varepsilon\in\mathcal E_M(\mathbb R)$, so
$u:=[(u_\varepsilon)_\varepsilon]\in\mathcal G(\mathbb R)$. One can easily see that for all $x\in\mathbb R$, $u_\varepsilon(x)=0$, whenever $\varepsilon$
is sufficiently large. Hence, $u(x)=0$ in $\widetilde{\mathbb R}$. But $u\neq 0$.

However, the following characterization holds (cf.\ \cite{MO1}):
\begin{theorem}\label{distingprop}
Let $f\in\mathcal G(\Omega)$. The following are equivalent:
\begin{enumerate}
\item $f=0$ in $\mathcal G(\Omega)$,
\item \label{comp} $(\forall\, x_c\,\in\widetilde{\Omega}_c)\,( \widetilde f(x_c)=0)$.
\end{enumerate}
\end{theorem}
Note that a similar statement holds in Egorov algebras (cf.\ \cite{Eyb4} and the final remark in \cite{MO1}). Furthermore,
Theorem \ref{distingprop} (\ref{comp}) can be weakened in the sense that only evaluation of $u$ at so-called near standard points
is necessary to ensure uniqueness in $\mathcal G(\Omega)$ (cf.\ (\cite{KKo1}, Proposition 4.2).

Finally, the composition $f\circ g$ of two generalized functions $f,\,g\in\mathcal G(\Omega)$ is defined, on the level of representatives as $f\circ g:=(f_\eps\circ g_\eps)_\eps+\mathcal N(\Omega)$. Note that composition in general requires $g$ to be
c-bounded (cf.\ \cite{KKo1}, section 2). However, if $h:=f\circ g$ is well-defined, then the evaluation map $\tilde h$
on $\widetilde{\Omega}_c$ agrees with the composition of the evaluation maps $\tilde f\circ\tilde g$. We shall use interchangeably
both forms of compositions. In this article, we shall restrict mostly to the cases $\Omega=\mathbb R^d$, $g(x)=x+h$, $g(x)=\lambda x$, $g(t)=b_c t$ with scaling parameter $\lambda\in\mathbb R^+$, $h\in\mathbb R^d$, $b_c\in\widetilde{\Omega}_c$ and $t\in\mathbb R$.

\subsection{Partial order on the ring of generalized numbers}
On the ring of generalized numbers a partial order can be defined as follows.
Let $\alpha\in\widetilde{\mathbb R}$. We say $\alpha\geq 0$, if there exists a representative
$(\alpha_\varepsilon)_\varepsilon$ such that for each $\varepsilon$ we have $\alpha_\varepsilon\geq 0$.
The partial ordering on the ring of generalized numbers therefore is given by the following
\[
\alpha\leq \beta\Leftrightarrow \beta-\alpha\geq 0.
\]
Note that contrary to the respective order on the real numbers, this order is not a total ordering. In addition to this order,
the so-called strict order on $\widetilde{\mathbb R}$ is used: Let $\beta\in \widetilde{\mathbb R}$. We say
$\beta$ is strictly positive (and write $\beta>>0$), if for each representative $(\beta_\varepsilon)_\varepsilon$
of $\beta$ there exists a number $m_0$ and an $\varepsilon_0\in(0,1]$ such that for each $\varepsilon<\varepsilon_0$
we have $\beta_\varepsilon>\varepsilon^{m_0}$. Note that the first two quantifiers can be interchanged. For further equivalent descriptions
of strict positivity, cf.\ \cite{algfound}. We say $\beta$ is strictly negative and write $\beta<<0$, if $-\beta>>0$. Furthermore, we shall write $a<<b$ whenever $a-b<<0$.
\subsection{Translation invariance in $\mathcal G$}
We shall use the following non-trivial fact (cf.\
\cite{Lobster}, \cite{Vernaeve} in the sequel:
 \begin{theorem}\label{lobster}
  Let $u\in\mathcal G(\mathbb R^d)$ such that
  \[
  (\forall\,h\in\mathbb R^d)(u(x+h)=u(x))
  \]
  holds in $\mathcal G(\mathbb R^d)$. Then $u$ is a constant.
   \end{theorem}

Theorem \ref{lobster} has found applications concerning characterizations of group invariants in algebras of generalized functions (\cite{KKo1,KKo2,Vernaeve}).

\section{Scaling invariance in $\mathcal G(\mathbb R^d)$}\label{mainsec}
To begin with we state the following:
\begin{lemma}\label{lemmaint}
Let $a,b\in\widetilde{\mathbb R}_c$ such that $a<<b$. Let
$f\in\mathcal G(\mathbb R)$. Then
\begin{enumerate}
\item \label{int1} the integral
\[
\int_a^b \vert f\vert\, dx:=\left(\int_{a_\varepsilon}^{b_\varepsilon}\vert f_\varepsilon\vert  \,dx\right)+\mathcal N(\mathbb R)
\]
is a well defined element of $\widetilde{\mathbb R}$. I.e., the
definition is independent of the choice of representatives
$(a_\varepsilon)_\varepsilon,\,
(b_\varepsilon)_\varepsilon,\,(f_\varepsilon)_\varepsilon$ of
$a,\,b,\,f$ respectively.
\item \label{int2} Furthermore, if $\int_a^b \vert f\vert=0$ in
$\widetilde{\mathbb R}$, then $f(x_c)=0$ for each generalized point
$x_c$ satisfying $a\leq x_c\leq b$.
\end{enumerate}
\end{lemma}
\begin{proof}
(\ref{int1}) Elementary.\\
(\ref{int2}) Assume $f(x_c)=c\neq 0$ for some $x_c$ with $a\leq x_c\leq b$. Let
$(x_{\varepsilon})_{\varepsilon}$ be a representative of $x_c$. As
$f(x_c)\neq 0$, there exists a zero sequence
$\varepsilon_k\rightarrow 0$ ($k\rightarrow\infty$) and an
$m_0\in{\mathbb N}$ such that $\vert
f_{\varepsilon_k}(x_{\varepsilon_k})\vert\geq 2\varepsilon_k^{m_0}$.
Applying the mean value theorem (and using the fact that the
derivative $(f'_\varepsilon)_\varepsilon$ is moderate), we therefore
have:
\begin{equation}
(\exists\,\rho\in\mathbb R^+)(\forall\, k\in\mathbb N)(\forall\,
y\in
[x_{\varepsilon_k}-\varepsilon_k^{\rho},\,x_{\varepsilon_k}+\varepsilon_k^{\rho}])(
\vert f_{\varepsilon_k}(y)\vert\geq\varepsilon_k^{m_0}).
\end{equation}
Choosing $\rho$ small enough, we can assume that $b_\varepsilon\ge
a_\varepsilon + \varepsilon^{\rho}$, $\forall \varepsilon$.
For sufficiently large $k$ we obtain
\begin{equation}\label{intcont}
\int_{a_{\varepsilon_k}}^{b_{\varepsilon_k}}\vert
f_{\varepsilon_k}(y)\vert \,dy\geq
\varepsilon_k^{m_0}(\varepsilon_k^{\rho})=\varepsilon_k^{\rho+m_0}.\end{equation}
Since $\left(\int_{a_\varepsilon}^{b_\varepsilon}\vert
f_{\varepsilon}(y)\vert\, dy\right)_{\varepsilon}$ is a
representative of $I(f)=\int_a^ b \vert f\vert dx$, inequality
(\ref{intcont})
contradicts our assumption $I(f)=0$ (the representative not being a negligible net).
\end{proof}
\begin{remark}\rm
 Note that if the condition $a<<b$ is weakened to $a\leq b$ and $a\neq b$ in (\ref{int2}) the conclusion of the above statement does not hold. To see this,
 we define $a,b\in\widetilde{\mathbb R}_c$ and $f\in\mathcal G(\mathbb R)$ on the level of representatives as follows:
 \[
 a_\varepsilon:=0,\; (\varepsilon\in(0,1]),\qquad b_\varepsilon:=\begin{cases} 1,\qquad \varepsilon=1/n \\ 0,\qquad\textit{otherwise}  \end{cases}
 \]
 and
 \[
 f_\varepsilon:=\begin{cases} 0,\qquad \varepsilon=1/n \\ 1,\qquad\textit{otherwise}  \end{cases}
 \]
 Then of course for each $\varepsilon>0$,
 \[
 \int_{a_\varepsilon}^{b_\varepsilon}\vert f_\varepsilon(x)\vert dx= 0,
 \]
 but $f\neq 0$ on the generalized interval $[a,b]$.
 \end{remark}
\begin{lemma}\label{prevlem}
Let $f\in\Gen(\R^+)$ and let $f(\lambda x)=f(x)$, $\forall \lambda\in\R^+$
(equality in $\Gen(\R^+)$). Then $f$ is a generalized constant, i.e., $f(x)=c$ holds
in $\Gen(\R^+)$, for some $c\in\widetilde{\C}$.
\end{lemma}
\begin{proof}
Write $g:=f\comp\exp$ (and check that this is well-defined in
$\Gen(\R)$). Then $g(x+h)=g(x)$, $\forall h\in\R$ (equality in $\Gen(\R)$). So
$g$ is a generalized constant in $\Gen(\R)$. So $f$ is a generalized constant in
$\Gen(\R^+)$.
\end{proof}
\begin{theorem}\label{ScTH}
Any generalized function $f$ in $\mathbb R^d$ which is invariant under standard scaling is constant.
\end{theorem}
\begin{proof}
We start by proving the statement in the one-dimensional case. As we shall see later,
the general case can be reduced to the latter.\\
{\it Case $d=1$}\\
According to Lemma \ref{prevlem}, $f$ is constant on $\mathbb R^+$ (that is $f=c$ in $\mathcal G(\mathbb R^+)$). By setting $\hat f(x):=f(-x)$ and applying
the preceding lemma again, we see that $f$ is also constant on $\mathbb R^-$. So $f$ is invariant under generalized scaling by $\lambda\in\widetilde{\mathbb R}^+_c=\widetilde{\mathbb R}_c\cap\widetilde{\mathbb R}^+$, 
\begin{equation}\label{non-infinitesimal scaling}
(\forall K\subset\subset\mathbb R^+)(\forall b\in\mathbb R)(\forall L\subset\subset\mathbb R^{+})
\left(\sup_{\lambda\in L}\sup_{x\in K}\vert f_\varepsilon(\lambda
x)-f_\varepsilon(x)\vert=O(\varepsilon^b),\,\varepsilon\rightarrow 0\right).
\end{equation}
In $\widetilde{\mathbb R}$,
\[\int_0^1\vert f(x)-f(1)\vert=\frac{1}{2}\int_0^2\vert f(x/2)-f(1)\vert.\]
By equation~(\ref{non-infinitesimal scaling}) for $K=\{1\}$,
$\int_1^2\vert f(x/2)-f(1)\vert=0$ in $\widetilde{\mathbb R}$.
Furthermore, since $f$ is scaling invariant, we have $\int_0^1\vert
f(x/2)-f(1)\vert= \int_0^1\vert f(x)-f(1)\vert$ in
$\widetilde{\mathbb R}$. We conclude that $\int_0^1\vert
f(x)-f(1)\vert=0$. By means of Lemma \ref{lemmaint} we conclude that
$f(0)=f(1)$. Similarly, $\int_{-1}^0 \vert f(x)-f(-1)\vert=0$.
Therefore, again by Lemma \ref{lemmaint} we have $f(0)=f(-1)$ and
$f(x)=f(0)$ on all generalized
points $x\in\widetilde{\mathbb R}_c$. \\
{\it Case $d>1$.}\\
Let $f\in\mathcal G(\mathbb R^d)$ invariant under positive (standard) scaling, that is,
$\forall\lambda\in\mathbb R$, $\lambda >0$ we have:
\[
f(\lambda x)=f(x).
\]
Fix a net $(a_\varepsilon)_\varepsilon$  such that $a_\varepsilon\in
L\subset\subset\mathbb R^d$ for all $\varepsilon>0$. Then the net
$(g_\varepsilon(t))_\varepsilon:=(f_\varepsilon(a_\varepsilon
t))_\varepsilon$ defines a generalized function
$g:=[(g_\varepsilon)_\varepsilon]\in\mathcal G(\mathbb R)$. Now the
scaling invariance for a fixed $\lambda$
 \[
(\forall\, K\subset\subset\mathbb R^d)(\forall\, b\in\mathbb
R)(\sup_{x\in K}\vert f_\varepsilon(\lambda x) -
 f_\varepsilon(x)\vert= O(\varepsilon^b),\ \text{ as } \varepsilon\rightarrow
 0)
\]
implies the scaling invariance for the same $\lambda$ of $g$
\[
(\forall K\subset\subset\mathbb R)(\forall\, b\in\mathbb
R)(\sup_{t\in K}\vert f_\varepsilon(\lambda a_\varepsilon t) -
f_\varepsilon(a_\varepsilon t)\vert= O(\varepsilon^b), \text{ as }
\varepsilon\rightarrow 0).
\]
So the one-dimensional case implies that $g$ is a generalized constant, that is,
\[
(\forall K\subset\subset\mathbb R)(\forall\, b\in\mathbb
R)(\sup_{t\in K}\vert f_\varepsilon(a_\varepsilon t) -
f_\varepsilon(0)\vert= O(\varepsilon^b), \text{ as
}\varepsilon\rightarrow 0).
\]
As $(a_\varepsilon)_\varepsilon$ is arbitrary, this implies that
\[
(\forall L\subset\subset\mathbb R^d)(\forall\, b\in\mathbb
R)(\sup_{x\in L}\vert f_\varepsilon(x) - f_\varepsilon(0)\vert=
O(\varepsilon^b), \text{ as }\varepsilon\to 0).
\]
Therefore $f$ is a generalized constant and we are done.
\end{proof}
\section{Homogeneity in the special algebra}\label{homsec}
We recall that a distribution $w\in\mathcal D'(\mathbb R^d) $
(resp.\ $\mathbb R^d\setminus\{0\}$) is homogeneous of degree
$\alpha\in\mathbb R$, if for all $\lambda\in\mathbb R^+$ and for all $\phi\in\mathcal
D(\mathbb R^d)$ (resp.\ $\mathcal D(\mathbb R^d\setminus\{0\})$) we
have
\[
\langle w(x),1/\lambda^d\phi(x/\lambda)\rangle=\lambda^\alpha\langle w,\phi\rangle.
\]

In this section we deal with two different notions of homogeneity in
$\mathcal G(\mathbb R^d)$. Our motivation is the following well
known fact in distribution theory (cf.\ Theorem 7.1.18 in
\cite{Hoermander1})
\begin{theorem}\label{homdist}
 Let $u\in\mathcal D'(\mathbb R^d\setminus\{0\})$ be homogeneous of degree $\alpha$. Then $u\in\mathcal S'(\mathbb R^d)$.
\end{theorem}
Before we go on to define an intrinsic notion of homogeneity, we
introduce tempered generalized functions. Let $\Omega$ be a
non-empty open subset of $\mathbb R^d$. The Colombeau algebra of
tempered generalized functions on $\Omega$ given by the quotient
\[
\mathcal G_\tau(\Omega):=\mathcal E_{M,\tau}(\Omega)/\mathcal
N_\tau(\Omega)
\]
where  the ring of tempered moderate nets of smooth functions is
given by
\[
\mathcal
E_{M,\tau}(\Omega):=\{(u_\varepsilon)_\varepsilon\in\mathcal
C^\infty(\Omega)^{(0,1]}\,|\,\exists\, N:\sup_{x\in\Omega}\vert
u_\varepsilon(x)\vert=O(\varepsilon^{-N}(1+\vert x\vert)^N))\}
\]
whereas ideal of tempered negligible functions is given by
\[
\mathcal N_\tau(\Omega):=\{(u_\varepsilon)_\varepsilon\in\mathcal
C^\infty(\Omega)^{(0,1]}\,|\,\exists\, N\,\forall\,
p:\sup_{x\in\Omega}\vert u_\varepsilon(x)\vert=O(\varepsilon^p
(1+\vert x\vert)^N))\}.
\]
The latter is an ideal in $\mathcal E_{M,\tau}(\Omega)$. Also, $\Gen_\tau(\Omega)$ can be linearly mapped into the special algebra
via the well-defined mapping
\[
K:\,\Gen_\tau(\Omega)\rightarrow \mathcal G(\Omega):\;\; ((K(u))_\varepsilon)_\varepsilon:=(u_\varepsilon)_\varepsilon+\mathcal N(\Omega)
\]
where $(u_\varepsilon)_\varepsilon$ is a representative of $u$.
This, however, is not an embedding (cf.\  \cite{C}, Proposition
4.1.6), since $K$ is not injective.
\subsection{Homogeneous generalized functions}
We start by introducing the notion of homogeneity in the special
Colombeau algebra.
\begin{definition}
A generalized function $u\in\mathcal G(\mathbb R^d)$ (resp.\
$\mathcal G(\mathbb R^d\setminus\{0\})$) is called homogeneous of
degree $\alpha$, if for each $\lambda\in\mathbb R^+$,
\begin{equation}
u(\lambda x)=\lambda^\alpha u(x)
\end{equation}
holds in $\mathcal G(\mathbb R^d)$ $($resp.\ $\mathcal G(\mathbb
R^d\setminus\{0\}))$.
\end{definition}
First we describe homogeneity of generalized functions outside the
origin of $\mathbb R^d$.
\begin{lemma}\label{keylemma}
Let $u\in\Gen(\R^+)$ and let $u(\lambda x)=\lambda^\alpha u(x)$,
$\forall \lambda\in\R^+$ and some $\alpha\in\R$ (equality in $\Gen(\R^+)$).
Then $u(x)=c x^\alpha$ in $\Gen(\R^+)$, for some $c\in\widetilde{\C}$.
\end{lemma}
\begin{proof}
Let $f(x):=\frac{u(x)}{x^\alpha}$. Then $f\in\Gen(\R^+)$. Moreover,
$f(\lambda x)=f(x)$, $\forall \lambda\in\R^+$ (equality in $\Gen(\R^+)$). By Lemma
\ref{prevlem}, $f(x)=c$ holds in $\Gen(\R^+)$, for some $c\in\widetilde{\C}$.
So $u(x) = c x^\alpha$ holds in $\Gen(\R^+)$.
\end{proof}
\begin{theorem}\label{maincharhom}
The following holds:
\begin{enumerate}
\item \label{441} Let $u\in\mathcal G(\mathbb
R^d\setminus\{0\})$ and let $u(\lambda x)=\lambda^\alpha u(x)$,
$\forall \lambda\in\R^+$ and some $\alpha\in\R$ (equality in
$\Gen(\R^d\setminus\{0\})$). Then $u(x)=
u\big(\frac{x}{\abs{x}}\big)\abs{x}^\alpha$ holds in
$\Gen(\R^d\setminus\{0\})$.
\item \label{442} Let $u\in\mathcal G(\mathbb R^d)$ be homogeneous of degree $\alpha$. Let $K$ be the
canonical map $\Gen_\tau(\R^d)\rightarrow\Gen(\R^d)$ as introduced above.
There exists $v\in\Gen_\tau(\R^d)$ such that $K(v)=u$. 
\end{enumerate}
\end{theorem}
\begin{proof}
(\ref{441}) Let $(u_\eps)_\eps$ be a representative of $u$. As for
each $\eps$,
$u_\eps\big(\frac{x}{\abs{x}}\big)\abs{x}^\alpha\in\Cnt\infty(\R^d\setminus\{0\})$
and $(u_\eps\big(\frac{x}{\abs{x}}\big)\abs{x}^\alpha)_\eps$
satisfies the moderateness-estimates, it represents an element of
$\Gen(\R^d\setminus\{0\})$. It remains to show that, for each
$K\csub\R^d\setminus\{0\}$ and $p\in\N$,
\begin{equation}\label{keyestimate}
\sup_{x\in
K}\abs{u_\eps(x)-u_\eps\left(\frac{x}{\abs{x}}\right)\abs{x}^\alpha}=O(\eps^p),\;
\eps\to 0.
\end{equation}
To this end, let $(a_\varepsilon)_\varepsilon$ represent a compactly supported point in $\mathbb R^d\setminus\{0\}$ and define $f\in\Gen(\R)$ on representatives by $f_\eps(t):=u_\eps(\frac{a_\eps}{\abs{a_\eps}} t)$.
Then $f(\lambda t)=\lambda^\alpha f(t)$ for all $\lambda\in\R^+$. So by Lemma \ref{keylemma},
$f(t)=ct^\alpha$ holds in $\Gen(\R^+)$, for some $c\in\widetilde{\mathbb C}$. Clearly,
$c=f(1)$. Set $(t_\varepsilon)_\varepsilon=(\abs{a_\eps})_\varepsilon$. Clearly, this represents a compactly
supported generalized constant
in $\R^+$. Therefore, for each $p\in\N$ we have
\[
\abs{u_\eps(a_\eps) -
u_\eps\left(\frac{a_\eps}{\abs{a_\eps}}\right)\abs{a_\eps}^\alpha} =
O(\eps^p),\, \eps\to 0.
\]
Since $(a_\varepsilon)_\varepsilon$ was arbitrary, we have shown (\ref{keyestimate}).

(\ref{442}) Take a representative $(u_\eps)_\eps$ of $u$. Choose a
cutoff function $\sigma\in\mathcal D(\mathbb R^d)$ which is
identically $1$ on $\Vert x\Vert\leq 1$. It is easily checked that
\[(\tilde u_\eps(x))_\eps:=\Big(u_\eps(x)\sigma(x)+
u_\eps\Big(\frac{x}{\abs{x}}\Big)\abs{x}^\alpha(1-\sigma(x))\Big)_\eps\]
is also a representative of $u$, and it coincides with
$(u_\varepsilon\big(\frac{x}{\abs{x}}\big)\abs{x}^\alpha)_\varepsilon$
on $\R^d\setminus \{x:\,\Vert x\Vert\leq 1\}$ as soon as $\eps$ is small enough. As
$u\in\Gen(\R^d)$, there exists for each $\beta\in\N^d$ an $N\in\N$
such that
\[\sup_{x\in\bar B(0,1)}\abs{\partial^\beta \tilde u_\eps(x)}=O(\eps^{-N}),\eps\to
0,\]
so in particular
\[\sup_{x\in\bar B(0,1)}(1+\abs{x})^{-N}\abs{\partial^\beta \tilde u_\eps(x)}=O(\eps^{-N}),\eps\to
0.\]
Further, for $\eps$ sufficiently small, and some $N\in\N$,
\begin{eqnarray}\nonumber
&&\sup_{x\in\R^d\setminus B(0,1)}(1+\abs{x})^{-N}\abs{\partial^\beta
\tilde u_\eps(x)} =\\\nonumber
&=&\sup_{x\in\R^d\setminus B(0,1)}(1+\abs{x})^{-N}\abs{\partial^\beta
\left(u_\eps\left(\frac{x}{\abs{x}}\right)\abs{x}^\alpha\right)}\le\\\nonumber
&\le&\sum_{\gamma\le\beta}\binom{\beta}{\gamma}
\underbrace
{\sup_{x\in\R^d\setminus
 B(0,1)}
 \abs{\partial^\gamma
 \left(u_\eps\left(\frac{x}{\abs{x}}
 \right)\right)}
 }_{O(\eps^{-N}),\eps\to 0}\times
 \\\nonumber
&\times&\underbrace{\sup_{x\in\R^d\setminus B(0,1)} (1+\abs{x})^{-N}\partial^{\beta-\gamma}
\abs{x}^\alpha}_{=O(1),\eps\to 0 \text{, if }\alpha\le N}.
\end{eqnarray}
So $v:=(\tilde u_\eps)_\eps + \Null_\tau(\R^d)$ satisfies the requirements.
\end{proof}\\\\
We note the following:
\begin{remark}\rm
\begin{enumerate}
\item In the one dimensional situation and for $\alpha$ not a negative integer, Theorem \ref{maincharhom} (\ref{441}) is analogue to the distributional setting, since every homogeneous distribution $w\in\mathcal D'(\mathbb R\setminus\{0\})$ can be written in the form $w(x)=C_1 x^\alpha_++C_2 x^\alpha_-$ (cf.\ \cite{Gelfand}, p.\ 87 or \cite{Hoermander1}, p.\ 75). Further they can be uniquely extended to the real line as homogeneous distributions of the same degree. However, this is not generally true in $\mathcal G$, cf.\ section Theorem \ref{superextension}, (\ref{extend2}), unless $\alpha\in\mathbb N_0$.
\item For a generalization of Theorem \ref{maincharhom} (\ref{442}), see Corollary \ref{mainex}.
\end{enumerate}

\end{remark}
What follows is a complete characterization of
 homogeneous generalized functions. We start with the following technical lemma:\\
  \begin{lemma}\label{technicality}
  Let $\widetilde{\mathbb C}[x_1,\dots,x_d]$ denote the ring of polynomials
 over the ring $\widetilde{\mathbb C}$ and let $k\in\mathbb N_0$. The following are equivalent:
  \begin{enumerate}
   \item \label{polytec1} $f$ lies in $\in\mathcal G(\mathbb R^d)\cap \widetilde{\mathbb C}[x_1,\dots,x_d]$ and is homogeneous of degree $k$, and
   \item \label{polytec2} For each $a_c\in\widetilde{\mathbb R}^d_c$ we have $g(t):=f(a_ct)\in\mathcal G(\mathbb R)\cap \widetilde{\mathbb C}[t]$ and
   is homogeneous of degree $k$.
  \end{enumerate}
 \end{lemma}
\begin{proof}
 Since the first implication is clear, we only need to prove (\ref{polytec2})$\Rightarrow$(\ref{polytec1}). For the sake of simplicity
we consider the case $k=1,\, d=2$. The general case is analogous. So
let $a_c\in\widetilde{\mathbb R}^d_c$, $g(t)=f(a_c t)\in\mathcal
G(\mathbb R)$. Since $g$ is a homogeneous polynomial of degree
 $k=1$, $g'$ is scaling invariant, therefore, $g$ is a constant. On the other hand, we have by the chain rule,
 \[
 g'(t)=\frac{\partial f}{\partial x_1}(a_ct)a_c^{(1)}+\frac{\partial f}{\partial x_2}(a_ct)a_c^{(2)}.
 \]
 This implies $g'(t)=g'(0)$ in $\mathcal G(\mathbb R)$, that is
 \[
 g'(t)=\frac{\partial f}{\partial x_1}(0)a_c^{(1)}+\frac{\partial f}{\partial x_2}(0)a_c^{(2)}.
 \]
 By integrating and homogeneity of $f$, we obtain
 \[
 f(a_ct)=g(t)=\frac{\partial f}{\partial x_1}(0)\{a_c^{(1)}t\}+\frac{\partial f}{\partial x_2}(0)\{a_c^{(2)}t\}.
 \]
 Setting $t=1$, we have
 \[f(a_c)=\frac{\partial f}{\partial x_1}(0)\{a_c^{(1)}\}+\frac{\partial f}{\partial x_2}(0)\{a_c^{(2)}\}\]
 and since $a_c$ was arbitrary, we have shown that $f=Ax_1+Bx_2$ with $A=\frac{\partial f}{\partial x_1}(0)$
 and $B=\frac{\partial f}{\partial x_2}(0)$ and we are done.
 \end{proof}\\\\
In the following theorem we establish that all homogeneous generalized functions are polynomials.
 \begin{theorem}\label{mainthhom}
Let $u\in\mathcal G(\mathbb R^d)$. Then,
\begin{enumerate}
 \item \label{strongconj1} if $u$ is a non-trivial homogeneous generalized function of degree $\alpha$, then $\alpha\in\mathbb N_0$.
 \item \label{strongconj2} As a consequence, the only non-trivial homogeneous generalized functions are homogeneous polynomials of degree $\alpha$ (in $d\geq 1$ variables).
\end{enumerate}
\end{theorem}
\begin{proof}
{\it Proof of (\ref{strongconj1}), case $d=1$.}\\
Consider $f_\eps(x):= u_\eps(x)/x^\alpha$. Then $f_\eps(x)\in\Cnt\infty(\R^+)$, $\forall\eps$.
Let $M\in\N$ fixed.
Let $(a_\eps)_\eps$ represent an element of $\widetilde\R_c$ with $a_\eps
\ge \eps^M$, $\forall\eps$. Let $g_\eps(x):= f_\eps(a_\eps x)$. As
$g_\eps\in\Cnt\infty(\R^+)$, $\forall\eps$, and $(g_\eps)_\eps$ satisfies the
moderateness-estimates on subcompacta of $\R^+$, it represents an element of $\Gen(\R^+)$.
Further, for $\lambda\in\R^+$ and $K\csub\R^+$,
\begin{eqnarray}\nonumber
\sup_{x\in K}\abs{g_\eps(\lambda x)-g_\eps(x)}
&\le&\sup_{x\in K}\frac{1}{\abs{\lambda a_\eps x}^\alpha}
\sup_{x\in K}\abs{u_\eps(\lambda a_\eps x) - \lambda^\alpha u_\eps(a_\eps x)}
=\\\nonumber&=& O(\eps^{-M\abs\alpha})\sup_{x\in L}\abs{u_\eps(\lambda x) - \lambda^\alpha u_\eps(x)},
\end{eqnarray}
for some $L\csub\R$, so $g(\lambda x)= g(x)$, $\forall \lambda\in\R^+$ (equality in $\Gen(\R^+)$),
so $g$ is a generalized constant in $\Gen(\R^+)$.
In particular, $\forall p\in\N$,
$\abs{g'_\eps(1)} = \abs{a_\eps}\abs{f'_\eps(a_\eps)} = O(\eps^p)$,
as $\eps \to 0$. So also, $\forall p\in\N$, $\abs{f'_\eps(a_\eps)}=O(\eps^p)$, as $\eps\to 0$.
As $(a_\eps)_\eps$ is arbitrary (with $a_\eps\ge\eps^M$), $\forall p\in\N$,
\[\sup_{x\in[\eps^M,M]}\abs{f_\eps'(x)}=O(\eps^p),\eps\to 0.\]
Then by Taylor's formula, $\forall p\in\N$,
\[
\sup_{x\in [\eps^M, M]} \abs{f_\eps(x)-f_\eps(1)}
\le \sup_{x\in [\eps^M, M]}\abs{x-1}\sup_{x\in [\eps^M, M]}\abs{f'_\eps(x)}
=O(\eps^p),\eps\to 0.
\]
Write $c_\eps:=f_\eps(1)$. Then also, $\forall p\in\N$,
\begin{equation}\label{alpha-in-N}
\sup_{x\in [\eps^M, M]} \abs{u_\eps(x)-c_\eps x^\alpha}
\le \sup_{x\in [\eps^M, M]}\abs{x^\alpha}\sup_{x\in [\eps^M, M]}\abs{f_\eps(x)-c_\eps}
=O(\eps^p),\eps\to 0.
\end{equation}
Now suppose that $\alpha < 0$. Then
\[
\sup_{x\in[0,1]}\abs{u_\eps(x)}
\ge\sup_{x\in[\eps^M,1]}\abs{u_\eps(x)}
\ge\abs{c_\eps}\sup_{x\in[\eps^M,1]}\abs{x^\alpha}
-\underbrace{\sup_{x\in[\eps^M,1]}\abs{u_\eps(x)-c_\eps x^\alpha}}
_{O(\eps^p),\eps\to 0, \forall p}
\]
As $M$ is arbitrary, this contradicts the moderateness of $u$, unless
$c_\eps=O(\eps^p)$, $\forall p\in\N$. So $\forall M$, $\forall p$,
\[
\sup_{x\in[\eps^M,M]}\abs{u_\eps(x)}=O(\eps^p), \eps\to 0.
\]
Then also, for $p\ge M+1$,
\begin{eqnarray}\nonumber
\sup_{x\in[0,M]}\abs{u_\eps(x)}
&\le& \sup_{x\in[0,M]}\abs{u_\eps(x)-u_\eps(x+\eps^p)}
+ \sup_{x\in[\eps^p,p]}\abs{u_\eps(x)}\\
\\\nonumber&\le& \eps^p\underbrace{\sup_{x\in[0,M+1]}\abs{u'_\eps(x)}}
_{O(\eps^{-N}),\eps\to 0}
+ \underbrace{\sup_{x\in[\eps^p,p]}\abs{u_\eps(x)}}_{O(\eps^l),\eps\to
0,\forall l}
\end{eqnarray}
for some $N\in\N$. Similarly, we have
$\sup_{x\in[-M,0]}\abs{u_\eps(x)}=O(\eps^p)$, $\forall p\in\N$. So $u=0$.\\
Now suppose that $\alpha>0$. Write $v_\eps(x):=u_\eps(x)-c_\eps x^\alpha$.
Then by Taylor's formula and by eqn.\ (\ref{alpha-in-N}),
\begin{eqnarray}\nonumber
\sup_{x\in[\eps^M,M]}\abs{v'_\eps(x)}
&\le&\eps^{-p}\sup_{x\in[\eps^M,M]}\abs{v_\eps(x+\eps^p)}
+\eps^{-p}\sup_{x\in[\eps^M,M]}\abs{v_\eps(x)}
+\\\nonumber&+&\frac{\eps^{p}}{2}
\underbrace{\sup_{x\in[\eps^M,M]}\abs{v''_\eps(x)}}_{O(\eps^{-N}),\eps\to 0},
\end{eqnarray}
for some $N\in\N$. So we also have, $\forall p\in\N$,
\[
\sup_{x\in[\eps^M,M]}\abs{u'_\eps(x)-\alpha c_\eps x^{\alpha-1}}=
O(\eps^p),\eps\to 0.
\]
If $0<\alpha<1$, this yields a similar contradiction with the moderateness of
$u'$, unless $u=0$. One proceeds inductively for $1<\alpha<2$, and so on.\newline\newline
{\it Proof of (\ref{strongconj1}), case $d>1$.}\\
If $u\ne 0$, there exists $(a_\eps)\in\widetilde{\R^d_c}$ such that
$(u_\eps(a_\eps))_\eps\ne 0$ as a generalized number. Then
$f_\eps(t):= u_\eps(a_\eps t)$ represents an element $f\in\Gen(\R)$
and $f(\lambda t)=\lambda^\alpha f(t)$, $\forall \lambda\in\R^+$
(equality in $\Gen(\R)$). Moreover, $f\ne 0$, since
$(f_\eps(1))_\eps\ne 0$ as a generalized number. So by the above,
$\alpha\in\N_0$.\newline\newline
{\it Proof of (\ref{strongconj2}).}\\
Let $u$ be an element of $\mathcal G(\mathbb R^d)$ that is homogeneous of degree $\alpha\in\mathbb R$. By
(\ref{strongconj1}), we know that $\alpha=k\in\mathbb N_0$. Furthermore, according to Lemma \ref{technicality}, it is sufficient to consider the case $d=1$. Differentiating $u$ $k$--times yields a scaling invariant $g\in\mathcal G(\mathbb R)$, therefore by Theorem \ref{ScTH}, $\exists c\in\widetilde{\mathbb R}$ such that
 $g=c$ in $\mathcal G(\mathbb R)$. Therefore $f\in \widetilde{\mathbb C}[x]$ is of degree $\deg f=k$, and by homogeneity, $f=cx^k$ in $\mathcal G(\mathbb R)$ and we are done.
\end{proof}
\subsection{The Euler Equation}
Aim of this section is to relate homogeneity in $\mathcal G(\mathbb R^d)$ to the Euler differential equation
\begin{equation}\label{eulereq}
\sum_{i=1}^d x_i\frac{\partial u}{\partial x_i}=\alpha u.
\end{equation}
We start by recalling the following fact in $\mathcal D'(\mathbb R^d)$ (cf.\ \cite{Gelfand}, pp.\ 286--287):
\begin{theorem}\label{eulerdistrichar}
Let $u\in \mathcal D'(\mathbb R^d)$, $\alpha\in\mathbb R$. The following are equivalent:
\begin{enumerate}
\item $u$ satisfies the Euler differential equation (\ref{eulereq}),
\item $u$ is homogeneous of degree $\alpha$.
\end{enumerate}
\end{theorem}
We shall prove the analogous statement in the setting of $\mathcal
G(\mathbb R^d)$. To start with, we note that we may confine
ourselves to the one dimensional case, that is $d=1$:
\begin{lemma}\label{prest}
Let $\Omega$ be either $\mathbb R^d$ or $\mathbb R^d\setminus\{0\}$. Let $u\in\mathcal G(\Omega)$, $\alpha\in\mathbb R$. The following are equivalent:
\begin{enumerate}
\item $u$ satisfies (\ref{eulereq}),
\item for each $b_c\in \widetilde{\Omega}_c$, $g(t):=u(b_c t)$ satisfies
\begin{equation}\label{eulerdim1}
t\frac{dg}{dt}(t)=\alpha g(t)
\end{equation}
in $\mathcal G(\mathbb R)$ (resp.\ in $\mathcal G(\mathbb R\setminus\{0\})$).
\end{enumerate}
\end{lemma}
\begin{proof}
Follows immediately by the chain rule and the fact that generalized
functions are determined by their values at compactly supported
generalized points.
\end{proof}\\\\
Next, we introduce a subset $J$ of $\mathcal G(\mathbb R)$ as
follows
\[
J:=\{g\in\mathcal G(\mathbb R)\,:\,(\forall u\in\mathcal G)\;(u=0\Leftrightarrow gu=0)\}.
\]
Hence, $J$ is precisely the set of all elements in $\mathcal
G(\mathbb R)$ which are not zero divisors, hence it is stable under
multiplication. We shall need the following fact to characterize all
solutions of the Euler Equation:
\begin{lemma}\label{pseudoinvertibility}
For each $k\in\mathbb N_0$, we have $x^k\in J$.
\end{lemma}
\begin{proof}
This is a simple extension of (\cite{KS1}, Example 2.3). See also the appendix.
\end{proof}
\begin{remark}\rm
\begin{enumerate}
\item
Note that the preceding statement is another confirmation of the
fact that the embedding $\iota$ of $\mathcal D'$ into $\mathcal G$
does not in general commute with multiplication of distributions
with smooth functions. To see this, consider a distribution
$w\in\mathcal D'(\mathbb R)$ with $\supp(w)={0}$. Then
$w=\sum_{j=0}^k c_j \delta^{(j)}$, where $c_j\in\mathbb C$. We have
$x^m w=0$, in $\mathcal D'(\mathbb R)$ for $m\geq k$ but $w\neq 0$.
On the one hand, $\iota(0)=\iota(x^kw)$, on the other hand
$\iota(x)\iota(w)\neq 0$ in $\mathcal G(\mathbb R)$. Hence
$\iota(xw)\neq\iota(x)\iota(w)=x\iota(w)$.
\item For a complete characterization of zero-divisors, see the appendix.
\end{enumerate}

\end{remark}
We are now able to prove our main statement:
\begin{theorem}
Let $\Omega$ be either $\mathbb R^d$ or $\mathbb R^d\setminus\{0\}$. 
Let $u\in\mathcal G(\Omega)$, $\alpha\in\mathbb R$. The following are equivalent:
\begin{enumerate}
\item \label{char3x} $u$ satisfies (\ref{eulereq}) in $\mathcal G(\Omega)$,
\item \label{char4x} $u$ is homogeneous of degree $\alpha$.
\end{enumerate}
\end{theorem}
\begin{proof}
{\bf Case $\Omega=\mathbb R^d$.}\\
According to Lemma \ref{prest}, we only need to show the claim for $d=1$.\\
(\ref{char4x})$\Rightarrow$(\ref{char3x}). By Theorem
\ref{mainthhom}, we have $\alpha$ is a non-negative integer, and $u$
is a monomial of degree $\alpha$ that is, $u=ax^\alpha$ with
$a=u(1)$, hence $u$ satisfies equation (\ref{eulereq}). \newline
(\ref{char3x})$\Rightarrow$(\ref{char4x}). Our strategy is to lead
$\alpha\in\mathbb R\setminus \mathbb N_0$ ad absurdum. For the sake
of convenience, we shall use the notation of the proof of Theorem
\ref{mainthhom}. Let $(u_\varepsilon)_\varepsilon$ be a
representative of $u$. Consider $f_\eps(x):= u_\eps(x)/x^\alpha$.
Then $f_\eps(x)\in\Cnt\infty(\R^+)$, $\forall\eps$. Let $M\in\N$ be
fixed. Let $(a_\eps)_\eps$ represent an element of $\widetilde\R_c$
with $a_\eps \ge \eps^M$, $\forall\eps$. Let $g_\eps(x):=
f_\eps(a_\eps x)$. As $g_\eps\in\Cnt\infty(\R^+)$, $\forall\eps$,
and $(g_\eps)_\eps$ satisfies the moderateness-estimates on
subcompacta of $\R^+$, it represents an element of $\Gen(\R^+)$.
Furthermore, $f$ solves $xf'=0$ in $\mathcal G(\mathbb R^+)$, hence
$f'=0$ in $\mathcal G(\mathbb R^+)$. Similarly, $g'=0$ in $\mathcal
G(\mathbb R^+)$. Hence $f$, $g$ are generalized constants on
$\mathbb R^+$. Now one can proceed as in the proof of Theorem
\ref{mainthhom} and obtains that $\alpha$ must be a non-negative
integer.\newline Let $\alpha=k=0$. Then the Euler equation reads
$xu'=0$, hence by Lemma \ref{pseudoinvertibility} we have $u'=0$,
hence $u$ must be a constant. If the order of homogeneity of $u$
equals $k=1$, then the differentiated Euler equation $xu'=u$ reads
$xu''=0$, hence, again by Lemma \ref{pseudoinvertibility} we
conclude that $u''=0$, hence $u=cx+d$. Inserting this solution into
$xu'=u$ yields $d=0$. As a consequence, $u$ is a monomial of degree
$1$, hence a homogeneous function of degree $1$. For degree of
homogeneity larger than $1$, one can proceed by induction. Hence $u$
is a monomial of degree $\alpha$. Thus  $u$ is homogeneous of degree
$\alpha$ and we are done with the first case.\\
{\bf Case $\Omega=\mathbb R^d\setminus\{0\}$.}\\
According to Lemma \ref{prest}, we may again confine ourselves to the case $d=1$. (\ref{char3x})$\Rightarrow$(\ref{char4x}): According to (Theorem 1.5.2 in \cite{Bible}), on $\mathbb R^+$ (resp.\ $\mathbb R^-$)
there exists a unique solution $v$ to the initial value problem (\ref{eulerdim1}) with initial data
$v(1)=u(1)$ (resp.\ $v(-1)=u(-1)$). Furthermore,
we know one solution, namely $v(x)=u(-1)(-x)^\alpha$ (resp.\ $v(x)=u(1)x^\alpha$). By uniqueness, $v=u$, hence $u$ is homogeneous of degree
$\alpha$.\\
(\ref{char4x})$\Rightarrow$(\ref{char3x}). Since $u$ is homogeneous of degree $\alpha$, by means of Theorem 
\ref{maincharhom} (\ref{441}), $u$ takes the form $u(x)=\begin{cases} u(-1)(-x)^\alpha, & x<0 \\ u(1)x^\alpha, & x>0\end{cases}$  . Hence $u$ satisfies
eq. (\ref{eulerdim1}) and we are done.
\end{proof}\\\\
In the following we make a comparison of the above with the distributional setting:
\begin{example}\rm
We consider homogeneity of degree $\alpha=-1$,
\begin{enumerate}\rm
\item In the distributional situation, we have by Theorem \ref{eulerdistrichar}, that
if  $u\in\mathcal D'(\mathbb R)$ is homogeneous of degree $-1$, then
\begin{equation}\label{example}
xu'+u=(xu)'=0.
\end{equation}
Hence, $xu=A$ in $\mathcal D'(\mathbb R)$, with $A$, an arbitrary constant. Now,
a particular solution of (\ref{example}) is $A\, 1/x$, with $1/x$, the principal value distribution. Hence the general solution
of (\ref{example}) is $u=A\,1/x+B\delta$ with $A,\,B\in\mathbb C$ (cf.\cite{Friedlander}, Theorem 2.7.1 and Exercise 2.3).
\item In the special algebra, $(xu)'=0$ implies $xu=c$ holds in $\mathcal G(\mathbb R)$ with $c$ a generalized constant.
However, then $c$ must be identically zero (because on the level of representatives, $xu=c$ violates moderateness of
representatives of $u$ whenever $c\neq 0$). Hence $xu=0$, therefore $u=0$ by Lemma \ref{pseudoinvertibility}.
\end{enumerate}
\end{example}
The preceding example suggests that the Euler equations should be solved with coupled calculus, meaning that apart from looking
for solutions of (\ref{eulerdim1}) one should solve this equation on the level of association, that is we solve for $g\in\mathcal G(\mathbb R)$
the equation
\begin{equation}\label{eulerdim1ass}
t\frac{dg}{dt}(t)\approx\alpha g(t).
\end{equation}
Consistency with the distributional setting is given via the following
\begin{theorem}
Let $w\in\mathcal D'(\mathbb R)$, $u=\iota(w)$ and $\alpha\in\mathbb R$. The following are equivalent:
\begin{enumerate}
\item $w$ solves $xw'=\alpha w$.
\item $u$ solves $xu'\approx \alpha u$.
\end{enumerate}
\end{theorem}
\begin{proof}
The proof follows from Proposition 1.2.70 (i) in \cite{Bible}.
\end{proof}
\subsection{Extension of homogeneous functions}
In this section extendability of homogeneous functions from the
pierced space to all of $\mathbb R^d$ is discussed. We shall say
$u\in\mathcal G(\mathbb R^d\setminus\{0\})$ is extendable, if there
exists $\hat u\in\mathcal G(\mathbb R^d)$ such that $\hat
u|_{(\mathbb R^d\setminus\{0\}} =u$ in $\mathcal G(\mathbb
R^d\setminus\{0\})$. Similarly, if $u\in\mathcal G(\mathbb
R^d\setminus\{0\})$ is homogeneous of degree $\alpha$, we say $u$ is
extendable as a homogeneous function, if there exists $\hat
u\in\mathcal G(\mathbb R^d)$ homogeneous of degree $\beta$ such that
$\hat u|_{(\mathbb R^d\setminus\{0\})} =u$ in $\mathcal G(\mathbb
R^d\setminus\{0\})$. Our first observation is
\begin{lemma}\label{lemmaex}
Let $u\in\mathcal G(\mathbb R^d\setminus\{0\})$, $u\neq 0$ be
homogeneous of degree $\alpha$. If $u$ is extendable as a
homogeneous function, then for the degree $\beta$ of homogeneity of
$\hat u$ we have $\alpha=\beta$.
\end{lemma}
\begin{proof}
Since $u\neq 0$, there exists $\tilde x_c\in (\widetilde{\mathbb
R^d\setminus\{0\}})_c$ such that $u(\tilde x_c)\neq 0$ in
$\widetilde{\mathbb R}$. Hence for $\lambda\in\mathbb R^+$, $0=u(\lambda\tilde
x_c)-u(\lambda\tilde x_c)=(\lambda^\alpha-\lambda^\beta)u(\tilde x_c)$. Hence
$\lambda^\alpha=\lambda^\beta$, which yields $\alpha=\beta$.
\end{proof}
\begin{theorem}\label{superextension}
We have the following:
\begin{enumerate}
\item \label{extend1}  Let $u\in\mathcal G(\mathbb R^d\setminus\{0\})$, $u\neq 0$ be homogeneous of degree $\alpha$. Then $u$ is extendable.
\item \label{extend2} Let $u\in\mathcal G(\mathbb R^d\setminus\{0\})$, $u\neq 0$ be homogeneous of degree $\alpha$. If $u$ is extendable as a homogeneous function,
then necessarily we have $\alpha\in\mathbb N_0$.
\end{enumerate}
\end{theorem}
\begin{proof}
{\it Proof of (\ref{extend1})}\\
Let $u\in\mathcal G(\mathbb R^d\setminus\{0\})$, $u\neq 0$ be
homogeneous of degree $\alpha$. By Theorem \ref{maincharhom}, $u$
can be written in the form
\[
u(x)=u\left(\frac{x}{\vert x\vert}\right)\vert x\vert ^\alpha.
\]
Let $\sigma\in\mathcal D(\mathbb R)$ be a cutoff function at $x=0$,
say $\sigma\equiv 1$ on $\Vert x\Vert\leq 1$ and $\sigma\equiv 0$ for $\Vert x\Vert\geq 2
$. Set $\rho_\varepsilon(x):=1-\sigma(\frac{\vert
x\vert }{\varepsilon})$. Then
$(\rho_\varepsilon)_\varepsilon\in\mathcal E_M(\mathbb R^d)$. Let
$(u_\varepsilon)_\varepsilon$ be a representative of $u$. We define
a net $\hat u_\varepsilon$ of smooth functions on $\mathbb R^d$ by
\[
\hat u_\varepsilon(x)=\rho_\varepsilon(x)u_\varepsilon\left(\frac{x}{\vert x\vert}\right)\vert x\vert ^\alpha.
\]
Clearly, $(\hat u_\varepsilon|_{\vert x\vert>0})_\varepsilon$ yields a representative of $u$ as well. It remains to prove that
$(\hat u_\varepsilon)_\varepsilon\in\mathcal E_M(\mathbb R^d)$. Fix $R>0$. First we estimate the zero derivative of $\hat u_\varepsilon$  on $B(0,R)=\{x:\Vert x\Vert\leq R\}$.
Since $\rho_\varepsilon(x)=0$ for $\vert x\vert <\varepsilon$, it suffices to estimate on $\varepsilon\leq \vert x\vert \leq R$. We have $\rho_\varepsilon(x)=O(1)$
on $\mathbb R^d$, $\vert x\vert^\alpha\leq\max\{\varepsilon^{-\vert\alpha\vert}, R^\alpha\}$ and for some $N$, $\vert \hat u_\varepsilon(\frac{x}{\vert x\vert})\vert=O(\varepsilon^{-N})$. Hence
\[
\sup_{\varepsilon\leq \vert x\vert \leq R}\vert \hat u_\varepsilon(x)\vert=O(\varepsilon^{-N-\vert \alpha\vert}),
\]
and we are done with the zero-order estimates. Bounds on higher order derivatives of $\hat u_\varepsilon$ are achieved similarly.\\
{\it Proof of (\ref{extend2})}\\This is a consequence of Lemma \ref{lemmaex} and Theorem \ref{mainthhom} (\ref{strongconj1}).
\end{proof}\\\\
\begin{remark}\rm
Note that for each  $\alpha\in\mathbb N_0$, there exist extendable homogeneous functions of degree $\alpha$ which are, however, not extendable as homogeneous functions.
Indeed, let $n$ be a non-negative integer.
Set $w(x):=x_+^n=x^n H(x)$, with $H$, the Heaviside function.
Consider the distribution $u(x_1,\dots,x_d):=w(x_1)\otimes
1\otimes\dots\otimes 1\in\mathcal D'(\mathbb R^d)$. Since $w$ is smooth away
from zero, $\iota(u)=\sigma(u)$ in $\vert x\vert>0$ because of the
sheaf-theoretic properties of the embedding $\iota$. Hence $\iota
(u)|_{\vert x\vert>0}\in\mathcal G(\mathbb R^d\setminus \{0\})$ is
homogeneous of degree $n$. Assume now, $\iota (u)|_{\vert x\vert>0}$
extends to a homogeneous function $\hat u$. Then, by Lemma
\ref{lemmaex}, the degree of homogeneity of $\hat u$ is $n$ as well.
But then it follows from Theorem \ref{mainthhom} that $\hat u$ is a
polynomial with generalized coefficients. This is impossible,
because $u$ is not a polynomial with generalized coefficients.
\end{remark}

As a consequence of Theorem \ref{maincharhom} and Theorem \ref{superextension} we have:
\begin{corollary}\label{mainex}
Let $u\in\mathcal G(\mathbb R^d\setminus\{0\})$ be homogeneous of degree $\alpha$. Let $K$ be the
canonical map $\Gen_\tau(\R^d)\rightarrow\Gen(\R^d)$ as introduced above.
There exists $v\in\Gen_\tau(\R^d)$ such that $K(v)|_{\mathbb R^d\setminus\{0\}}=u$. 
\end{corollary}
\begin{proof}
Take $u$ and extend it by means of Theorem \ref{superextension} to $\hat u\in\mathcal G(\mathbb R^d)$.
The rest follows from the fact that in the proof of Theorem \ref{maincharhom} (\ref{442}), only homogeneity of
$\hat u$ in $\mathbb R^d\setminus\{0\}$ is used, which is the case according to our assumptions.
\end{proof}
\subsection{Weak homogeneity in generalized function algebras}
Aim of this section is to introduce a weaker concept of homogeneity
in $\mathcal G(\mathbb R^d)$ than the one discussed so far. Then we
shall compare these two notions of homogeneity and show that weak
homogeneity is consistent with the distributional notion of
homogeneity on the level of embedded distributions.
\begin{definition}
We call $u\in\mathcal G(\mathbb R^d)$ weakly homogeneous of degree $\alpha$, if
for each $\lambda\in\mathbb R^+$ and for each $\phi\in \mathcal D(\mathbb R^d\setminus\{0\})$ the
identity
$$ \int_{\mathbb R^d} u(\lambda x) \phi (x) dx=\int_{\mathbb R^d} \lambda^\alpha u(x)\phi(x) dx
$$
holds in $\widetilde{\mathbb R}$.
\end{definition}
\begin{remark}\rm\label{remhom1}
Clearly homogeneity in $\mathcal G(\mathbb R^d)$ implies weak
homogeneity.
\end{remark}

\begin{proposition}\label{distshadow}
Suppose we are given a generalized function $u\in\mathcal G(\mathbb R^d)$ such that for some $w\in\mathcal D'$, $\iota(w)=u$. Let $\alpha\in\mathbb R$. Then the following are equivalent:
\begin{enumerate}
 \item \label{distshadowitem1} $u$ is weakly homogeneous in $\mathcal G(\mathbb R^d)$ of degree $\alpha$,
 \item \label{distshadowitem2} $w$ is homogeneous in $\mathcal D'(\mathbb R^d)$ of degree $\alpha$.
\end{enumerate}
 \end{proposition}
 \begin{proof}
Since $\iota(w)=u$, the assertion is a consequence of (\cite{Bible}, Theorem 1.2.63): The identity
\[
(\forall\varphi\in\mathcal D(\mathbb R^d))(\int_{\mathbb R^d} u(x)\varphi(x) \,dx=\langle w,\varphi\rangle)
\]
holds in $\widetilde{\mathbb R}$.
\end{proof}
\begin{example}
 As a non-trivial example for the situation described in Proposition \ref{distshadow} let us consider the $\delta$ distribution. Let $\rho\in \mathcal S$
 be a mollifier allowing for an embedding $\iota:\, \mathcal D'(\mathbb R^d)\hookrightarrow \mathcal G(\mathbb R^d)$. Define
 $u:=(\rho_\varepsilon)_\varepsilon+\mathcal N(\mathbb R^d)$. Clearly, $\iota(\delta)=u$. Furthermore $\delta$ is homogeneous of degree
 $\alpha=-d$, therefore we have
\[
\int \rho_\varepsilon(\lambda x)\phi(x) dx+\mathcal N(\mathbb R^d)=
\lambda^{-d}\int \rho_\varepsilon(x)\phi(x) dx + \mathcal N(\mathbb R^d).
\]
\end{example}
It should be noted, that the converse of Remark \ref{remhom1} is not
true:
\begin{proposition}\label{moreweaksolutions}
The following holds:
\begin{enumerate}
\item \label{befc1} Weak homogeneity in $\mathcal G(\mathbb R^d)$ does not imply homogeneity in $\mathcal G(\mathbb R^d)$. Moreover, there exists $w\in\mathcal D'(\mathbb R^d)$ homogeneous of degree $\alpha$, but $\iota(w)$ is not homogeneous in $\mathcal G(\mathbb R^d)$.
\item \label{befc2} There exist generalized functions $u$, weakly homogeneous of degree $\alpha$, however not associated to a distribution.
\end{enumerate}
 \end{proposition}
\begin{proof} We start with the proof of (\ref{befc1}). For the sake of simplicity we consider first
$d=1,\alpha=0$. Take the Heaviside function $H\in\mathcal D'(\mathbb
R)$. Clearly $H$ is scaling invariant, that is $H$ is weakly
homogeneous of degree $\alpha=0$, and so is $\iota(H)$ according to Proposition \ref{distshadow}. 
Assume now, $\iota(H)$ is homogeneous of degree $\alpha=0$. This means that $\iota(H)$ is scaling invariant. By
Theorem \ref{ScTH}, $\iota(H)$ is a constant in $\mathcal G(\mathbb R)$. However, this is impossible, since
$\frac{d}{dx}\iota(H)=\iota\left(\frac{d}{dx}H\right)=\iota(\delta)\neq 0$ in $\mathcal G(\mathbb R)$. Therefore, $H$ is not homogeneous of the same degree.\\
The case $d\geq 1,\alpha=0$ can be shown by taking
$w=H\otimes\dots\otimes H$. Clearly, $w$ is scaling invariant, since
for $\varphi_i\in\mathcal D(\mathbb R)$, $1\leq i\leq d$, we have
$\langle
w,\frac{1}{\lambda^d}\varphi_1(x_1/\lambda)\otimes\dots\otimes\varphi_d(x_d/\lambda)\rangle=\prod_{i=1}^d\int_0^\infty\frac{1}{\lambda}\varphi_i(x/\lambda)dx=\prod_{i=1}^d\langle
H,\varphi_i\rangle$. However,
$\partial_1\iota(w)=\iota(\delta\otimes H\otimes\dots\otimes H)\neq
0$. Therefore, $\iota(w)$ cannot be scaling invariant in $\mathcal
G(\mathbb R^d)$ and we are done with $\alpha=0$. For degree of
homogeneity $\alpha\neq 0$, one proceeds similarly.\\
Proof of (\ref{befc2}). Take the net of real numbers
$(a_\varepsilon)_\varepsilon$ defined by $a_{1/n}=1$ ($n\in\mathbb
N$) and $a_\varepsilon=0$ otherwise. Then for
$a:=[(a_\varepsilon)_\varepsilon]$, $u:=a$ is a constant in $\mathcal
G(\mathbb R^d)$, hence it is homogeneous of degree $\alpha=0$. As
a consequence, $u$ is weakly homogeneous of the same degree, but $u$
is not associated to a distribution (in particular, it is not
associated to $0$ or $1$).
\end{proof}\\
We  conjecture the following analogue to Theorem \ref{homdist}:
\begin{conjecture}
Let $u\in\mathcal G(\mathbb R^d)$ be weakly homogeneous of degree
$\alpha$. Then there is a representative
$(u_\varepsilon)_\varepsilon$ of $u$ such that for any
$\phi\in\mathcal S(\mathbb R^d)$, the integral
\[
\left(\int_{\mathbb R^d} u_\varepsilon(x)\phi(x)dx\right)
\]
 yields a moderate net.
\end{conjecture}
\subsection{Associative homogeneity in the special algebra} The
notion ``weak homogeneity'' in $\mathcal G(\mathbb R^d)$ cannot in
general cope with generalized functions which admit a distributional
shadow which is homogeneous in $\mathcal D'(\mathbb R^d)$. In view
of Proposition \ref{distshadow} this means that weak homogeneity in
$\mathcal G$ is not consistent with homogeneity in $\mathcal D'$ on
the level of association. As an example, let $u\in\mathcal G(\mathbb
R)$ be defined by the class of $(u_\varepsilon)_\varepsilon$ where
$u_\varepsilon:=x^2+\varepsilon$. On the one hand, Proposition
\ref{distshadow} is not applicable, because $u\neq \iota(w)$, that
is, $u$ cannot be an embedded distribution. It is further clear that
$u$ is not weakly homogeneous of degree $2$: to see this, let
$\varphi\in\mathcal D(\mathbb R)$ with $\int \varphi dx=1$. Then
\[
\int u_\varepsilon\varphi dx=\int x^2\varphi dx+\varepsilon.
\]
In particular, for all $\lambda\neq 1$ we have
\[
\int (u_\varepsilon(\lambda x)-\lambda^2u_\varepsilon(x))\varphi
dx=\varepsilon(1-\lambda^2)\neq O(\varepsilon^2),
\]
and we have shown that $u$ is not weakly homogeneous of degree $2$.\\On the other hand, the distributional shadow $x^2$
clearly is homogeneous of degree $2$. This suggests the following intrinsic notion of an even weaker homogeneity in $\mathcal G$:
\begin{definition}
 Let $u\in\mathcal G(\mathbb R^d)$. We call $u$  associatively homogeneous
 of degree $\alpha$, if for all $\varphi\in\mathcal D(\mathbb R^d)$ and all $\lambda\in\mathbb R^+$, we have
 \[
  \int (u_\varepsilon(\lambda x)-\lambda^\alpha u_\varepsilon(x))\varphi(x)\,dx\rightarrow 0,\;(\varepsilon\rightarrow 0).
  \]
Note that this limit is independent of the choice of representative $(u_\varepsilon)_\varepsilon$ of $u$.
\end{definition}
The following is an analogue of Proposition \ref{distshadow} in the
context of this subsection. We skip the proof.
\begin{proposition}
 Let $u\in\mathcal G(\mathbb R^d)$ be given, with a distributional shadow $w\in\mathcal D'(\mathbb R^d)$. The following are
 equivalent:
 \begin{enumerate}
  \item $u$ is associatively homogeneous of degree $\alpha$.
  \item $w$ is homogeneous of degree $\alpha$.
 \end{enumerate}
\end{proposition}

\section*{Appendix. Zero-divisors in the special Colombeau algebra}
Let $\Omega$ be an open subset of $\mathbb R^d$. It is well known that in $\mathcal C(\Omega)$
zero divisors are precisely such functions which vanish on some non-empty open subset of $\Omega$.
In the special algebra $\mathcal G(\mathbb R)$ we have the following analog
\begin{theorem}
Let $f\in \mathcal G(\mathbb R)$, $f\neq 0$. The following are equivalent:
\begin{enumerate}
\item \label{one} $f$ is a zero divisor.
\item \label{two} Let $(f_\varepsilon)_\varepsilon$ be a representative of $f$. Then we have
\begin{eqnarray}\nonumber
&(\exists K\subset\subset \mathbb R)(\exists (x_\varepsilon)_\varepsilon\subset K^{(0,1]})(\exists\rho\geq 0)(\exists (\varepsilon_n)_n,\,\varepsilon_n\rightarrow 0)\\\nonumber&\qquad(\forall k\geq 0)(\forall x\in [x_{\varepsilon_k}-\varepsilon_k^\rho,x_{\varepsilon_k}+\varepsilon_k^\rho],\;\vert f_{\varepsilon_k}(x)\vert<\varepsilon_k^k).
\end{eqnarray}
\end{enumerate}
\end{theorem}
\begin{proof}
(\ref{two})$\Rightarrow$(\ref{one}). Take $\sigma\in C^\infty(\mathbb R)\setminus\{0\}$, such that $\supp \sigma\subseteq[-1,1]$. We define the net
$(\sigma_\varepsilon)_\varepsilon$ by
\[
\sigma_\varepsilon(x):=\begin{cases}\sigma(\frac{x-x_\varepsilon}{\varepsilon^\rho}),\;\varepsilon=\varepsilon_k\\0,\;\textit{otherwise}\end{cases}.
\]
We have $(\sigma_\varepsilon)_\varepsilon\in\mathcal E_M(\mathbb R)$ and $\sigma_\varepsilon=0$ whenever $\varepsilon\neq\varepsilon_k$
and whenever $\varepsilon=\varepsilon_k$ and $x\notin [x_{\varepsilon_k}-\varepsilon_k^\rho,x_{\varepsilon_k}+\varepsilon_k^\rho]$. Hence, with $g:=[(\sigma_\varepsilon)_\varepsilon]$
we have $fg=0$. But $g\neq 0$. Hence $f$ is a zero divisor in $\mathcal G(\mathbb R)$.\\
(\ref{one})$\Rightarrow$(\ref{two}). Let $0\neq g\in\mathcal G(\mathbb R)$ be given such that $fg=0$. Then by Theorem \ref{distingprop} there exists $x_c\in\widetilde{\mathbb R}_c$
such that $g(x_c)\neq 0$. In terms of representatives this means that 
\[
(\exists\varepsilon_k\rightarrow 0)(\exists m\geq 0)(\forall k\geq 0, \vert g_{\varepsilon_k}(x_{\varepsilon_k})\vert\geq 2\varepsilon_k^m).
\]
By the mean value theorem and moderateness of the first derivative, one even gets
\[
(\exists\varepsilon_k\rightarrow 0)(\exists m\geq 0)(\exists\rho\geq 0)(\exists k_0)(\forall k\geq k_0)(\forall x\in [x_{\varepsilon_k}-\varepsilon_k^\rho,x_{\varepsilon_k}+\varepsilon_k^\rho],\;\vert g_{\varepsilon_k}(x)\vert\geq \varepsilon_k^m).
\]
But $fg=0$, hence $f$ must satisfy
\[
(\forall l)(\forall x\in [x_{\varepsilon_k}-\varepsilon_k^\rho,x_{\varepsilon_k}+\varepsilon_k^\rho],\;\vert f_{\varepsilon_k}(x)\vert=O( \varepsilon_k^l),\;(k\rightarrow\infty)).
\]
Hence, by taking a subsequence of $\varepsilon_k$ the assertion follows.
\end{proof}
\section*{Acknowledgement}
The authors acknowledge funding by the Austrian Science Fund (FWF) through the research grants
P16742-N04, Y237-N13. Further, Hans Vernaeve acknowledges support by the Lise Meitner project grant M949-N18.

\end{document}